\documentclass{article}
\usepackage{amsmath,amsthm,amsopn,amstext,amscd,amsfonts,amssymb}
\usepackage{natbib,color}

\def \G {\mathop{\rm G}\nolimits}
\def \Q {\mathop{\rm Q}\nolimits}
\def \tr {\mathop{\rm tr}\nolimits}
\def \re {\mathop{\rm Re}\nolimits}
\def \im {\mathop{\rm Im}\nolimits}

\def \etr {\mathop{\rm etr}\nolimits}
\def \diag {\mathop{\rm diag}\nolimits}

\renewenvironment{abstract}
                 {\vspace{6pt}
                  \begin{center}
                  \begin{minipage}{5in}
                  \centerline{\textbf{Abstract}}
                  \noindent\ignorespaces
                 }
                 {\end{minipage}\end{center}}
\newtheorem{theorem}{\textbf{Theorem}}[section]

\newtheorem{lemma}{\textbf{Lemma}}[section]

\theoremstyle{definition}
\newtheorem{definition}{\textbf{Definition}}[section]

\newtheorem{remark}{\textbf{Remark}}[section]

\title{\Large \textbf{Distributions on symmetric cones I: Riesz distribution}}
\author{
  \textbf{Jos\'e A. D\'{\i}az-Garc\'{\i}a} \thanks{Corresponding author\newline
   {\bf Key words.}  Wishart distribution; Riesz distribution, gamma distribution, real, complex, quaternion
    and octonion random matrices.\newline
    2000 Mathematical Subject Classification. Primary 60E05, 62E15; secondary
    15A52}\\
  {\normalsize Department of Statistics and Computation} \\
  {\normalsize 25350 Buenavista, Saltillo, Coahuila, Mexico} \\
  {\normalsize E-mail: jadiaz@uaaan.mx} \\
}
\date{}
\begin{document}
\maketitle

\begin{abstract}
The Riesz distribution for real normed division algebras is derived in this work. Then two
versions of these distributions are proposed and some of their properties are studied.
\end{abstract}

\section{Introduction}\label{sec1}

In the context of the Statistics and related fields, the analysis on symmetric cones has been
developed in depth through the last seven decades. As usual, the studies were focused on the
real case, later the complex case appeared and in the recent years, the quaternionic and
octonionic numbers were included in the analysis. This analysis has been notorious in the
development of certain statistical areas as random matrix distributions, but also in some
mathematical topics including special polynomials matrix arguments as zonal, Jack, Davis'
invariants, Laguerre, Hermite, Hayakawa; and the associated hypergeometric type series. These
mathematical tools have been implemented for example in the characterisation of non-central
distributions, statistical theory of shape, etc.  Among many authors, we can highlight the
works of \citet{h:55}, \citet{j:61}, \citet{j:64}, \citet{c:63}, \citet{da:80}, \citet{m:82},
\citet{dg:11}, and the references therein. Of course, the interest in the analysis of symmetric
cones does not belong exclusively to statistics; the harmonic analysis, for example, is an
special case of the real symmetric cone and this cone play a fundamental role in number theory.
Another particular case that has been studied in detail, specially from the point of view of
the wave equation, is the Lorentz cone, \citet{fk:94} and the references therein.

Although during the 60's and 70's were obtained important results in the general analysis of
the symmetric cones,  the  past 20 years have reached a substantial development. Essentially,
these advances have been archived through two approaches based on the theory of Jordan algebra
and the real normed division algebras. A basic source of the general theory of symmetric cones
under Jordan algebras can be found in  \citet{fk:94}; and specifically, some works in the
context of distribution theory in symmetric cones based on Jordan algebras are provided in
\citet{m:94}, \citet{cl:96}, \citet{hl:01}, and \citet{hlz:05}, and the references therein. In
the field of spherical functions (Jack polynomials including James' zonal polynomials) we can
mention the work of \citet{s:97}. Parallel results on distribution theory based on real normed
division algebras have been also developed in statistics and random matrix theory, see \citet
{f:05}, \citet{dggj:11}, \citet{dggj:12a}, among others. A detailed study of functions with
matrix arguments on symmetric cones, such as Jack and Davis' invariants polynomials and
hypergeometric types series, among others, appear in \citet{gr:87}, \citet{dggj:12b} and
\citet{dggj:09b}, for real normed division algebras.

From the standpoint of multivariate statistics, the main focus of statistical studies has been
developed for real random matrices. However, many works have proposed extending some of these
tests to the case of complex random matrices, see for example \citet{w:56}, \citet{g:63},
\citet{j:64}, and \citet{mdm:06}, among many others. Studies examining the case of quaternionic
random matrices are much less common, see \citet{b:00}, \citet{lx:09} and \citet{dggj:13}. From
a applied point of view, the relevance of \emph{the octonions} remains unclear. An excellent
review of the history, construction and many other properties of octonions is given in
\citet{b:02}, where it is stated that:
\begin{center}
\begin{minipage}[t]{4in}
\begin{sl}
``Their relevance to geometry was quite obscure until 1925, when \'Elie Cartan described
`triality' -- the symmetry between vector and spinors in 8-dimensional Euclidian space. Their
potential relevance to physics was noticed in a 1934 paper by Jordan, von Neumann and Wigner on
the foundations of quantum mechanics...Work along these lines continued quite slowly until the
1980s, when it was realized that the octionions explain some curious features of string
theory... \textbf{However, there is still no \emph{proof} that the octonions are useful for
understanding the real world}. We can only hope that eventually this question will be settled
one way or another."
\end{sl}
\end{minipage}
\end{center}

For the sake of completeness, in the present article the case of octonions is considered, but
the veracity of the results obtained for this case can only be conjectured. Nonetheless,
\citet[Section 1.4.5, pp. 22-24]{f:05} it is proved that the bi-dimensional density function of
the eigenvalue, for a Gaussian ensemble of a $2 \times 2$ octonionic matrix, is obtained from
the general joint density function of the eigenvalues for the Gaussian ensemble, assuming $m =
2$ and $\beta = 8$, see Section \ref{sec2}.

Under the approach based on the theory of Jordan algebras, a family of distributions on
symmetric cones, termed the Riesz distributions, was first introduced by \citet{hl:01} under
the name of Riesz natural exponential family (Riesz NEF); it was based on a special case of the
so-called Riesz measure from \citet[p.137]{fk:94}, going back to \citet{r:49}. This Riesz
distribution generalises the matrix multivariate gamma and Wishart distributions, containing
them as particular cases.

This article introduces the Riesz distribution for real normed division algebras. Section
\ref{sec2} revises some definitions and notations on real normed division algebras, also, two
definitions of the generalised gamma funcion on symmetric cones are given; in addition, several
Jacobians with respect to Lebesgue measure for real normed division algebras are proposed.
Finally, two Riesz distributions are obtained and under both definitions, the corresponding
characteristic function and the eigenvalue distribution are derived in section \ref{sec3}.

\section{Preliminary results}\label{sec2}

A detailed discussion of real normed division algebras may be found in \citet{b:02} and
\citet{gr:87}. For convenience, we shall introduce some notation, although in general we adhere
to standard notation forms.

For our purposes, a \textbf{vector space} is always a finite-dimensional module over the field
of real numbers. An \textbf{algebra} $\mathfrak{F}$ is a vector space that is equipped with a
bilinear map $m: \mathfrak{F} \times \mathfrak{F} \rightarrow \mathfrak{F}$ termed
\emph{multiplication} and a nonzero element $1 \in \mathfrak{F}$ termed the \emph{unit} such
that $m(1,a) = m(a,1) = a$. As usual, we abbreviate $m(a,b) = ab$ as $ab$. We do not assume
$\mathfrak{F}$ associative. Given an algebra, we freely think of real numbers as elements of
this algebra via the map $\omega \mapsto \omega 1$.

An algebra $\mathfrak{F}$ is a \textbf{division algebra} if given $a, b \in \mathfrak{F}$ with
$ab=0$, then either $a=0$ or $b=0$. Equivalently, $\mathfrak{F}$ is a division algebra if the
operation of left and right multiplications by any nonzero element is invertible. A
\textbf{normed division algebra} is an algebra $\mathfrak{F}$ that is also a normed vector
space with $||ab|| = ||a||||b||$. This implies that $\mathfrak{F}$ is a division algebra and
that $||1|| = 1$.

There are exactly four normed division algebras: real numbers ($\Re$), complex numbers
($\mathfrak{C}$), quaternions ($\mathfrak{H}$) and octonions ($\mathfrak{O}$), see
\citet{b:02}. We take into account that $\Re$, $\mathfrak{C}$, $\mathfrak{H}$ and
$\mathfrak{O}$ are the only normed division algebras; moreover, they are the only alternative
division algebras, and all division algebras have a real dimension of $1, 2, 4$ or $8$, which
is denoted by $\beta$, see \citet[Theorems 1, 2 and 3]{b:02}. In other branches of mathematics,
the parameters $\alpha = 2/\beta$ and $t = \beta/4$ are used, see \citet{s:97} and
\citet{k:84}, respectively.

Let ${\mathcal L}^{\beta}_{m,n}$ be the set of all $n \times m$ matrices of rank $m \leq n$
over $\mathfrak{F}$ with $m$ distinct positive singular values, where $\mathfrak{F}$ denotes a
\emph{real finite-dimensional normed division algebra}. Let $\mathfrak{F}^{m \times n}$ be the
set of all $m \times n$ matrices over $\mathfrak{F}$. The dimension of $\mathfrak{F}^{m \times
n}$ over $\Re$ is $\beta mn$. Let $\mathbf{A} \in \mathfrak{F}^{n \times m}$, then
$\mathbf{A}^{*} = \overline{\mathbf{A}}^{T}$ denotes the usual conjugate transpose. The set of
matrices $\mathbf{H}_{1} \in \mathfrak{F}^{n \times m}$ such that
$\mathbf{H}_{1}^{*}\mathbf{H}_{1} = \mathbf{I}_{m}$ is a manifold denoted ${\mathcal
V}_{m,n}^{\beta}$, termed the \emph{Stiefel manifold} ($\mathbf{H}_{1}$ are also known as
\emph{semi-orthogonal} ($\beta = 1$), \emph{semi-unitary} ($\beta = 2$), \emph{semi-symplectic}
($\beta = 4$) and \emph{semi-exceptional type} ($\beta = 8$) matrices, see \citet{dm:99}). The
dimension of $\mathcal{V}_{m,n}^{\beta}$ over $\Re$ is $[\beta mn - m(m-1)\beta/2 -m]$. In
particular, ${\mathcal V}_{m,m}^{\beta}$ with dimension over $\Re$, $[m(m+1)\beta/2 - m]$, is
the maximal compact subgroup $\mathfrak{U}^{\beta}(m)$ of ${\mathcal L}^{\beta}_{m,m}$ and
consists of all matrices $\mathbf{H} \in \mathfrak{F}^{m \times m}$ such that
$\mathbf{H}^{*}\mathbf{H} = \mathbf{I}_{m}$. Therefore, $\mathfrak{U}^{\beta}(m)$ is the
\emph{real orthogonal group} $\mathcal{O}(m)$ ($\beta = 1$), the \emph{unitary group}
$\mathcal{U}(m)$ ($\beta = 2$), the \emph{compact symplectic group} $\mathcal{S}p(m)$ ($\beta =
4$) or \emph{exceptional type matrices} $\mathcal{O}o(m)$ ($\beta = 8$), for $\mathfrak{F} =
\Re$, $\mathfrak{C}$, $\mathfrak{H}$ or $\mathfrak{O}$, respectively. Denote by ${\mathfrak
S}_{m}^{\beta}$ the real vector space of all $\mathbf{S} \in \mathfrak{F}^{m \times m}$ such
that $\mathbf{S} = \mathbf{S}^{*}$. Let $\mathfrak{P}_{m}^{\beta}$ be the \emph{cone of
positive definite matrices} $\mathbf{S} \in \mathfrak{F}^{m \times m}$. Thus,
$\mathfrak{P}_{m}^{\beta}$ consist of all matrices $\mathbf{S} = \mathbf{X}^{*}\mathbf{X}$,
with $\mathbf{X} \in \mathfrak{L}^{\beta}_{m,n}$; then $\mathfrak{P}_{m}^{\beta}$ is an open
subset of ${\mathfrak S}_{m}^{\beta}$. Over $\Re$, ${\mathfrak S}_{m}^{\beta}$ consist of
\emph{symmetric} matrices; over $\mathfrak{C}$, \emph{Hermitian} matrices; over $\mathfrak{H}$,
\emph{quaternionic Hermitian} matrices (also termed \emph{self-dual matrices}) and over
$\mathfrak{O}$, \emph{octonionic Hermitian} matrices. Generically, the elements of
$\mathfrak{S}_{m}^{\beta}$ are termed as \textbf{Hermitian matrices}, irrespective of the
nature of $\mathfrak{F}$. The dimension of $\mathfrak{S}_{m}^{\beta}$ over $\Re$ is
$[m(m-1)\beta+2m]/2$. Let $\mathfrak{D}_{m}^{\beta}$ be the \emph{diagonal subgroup} of
$\mathcal{L}_{m,m}^{\beta}$ consisting of all $\mathbf{D} \in \mathfrak{F}^{m \times m}$,
$\mathbf{D} = \diag(d_{1}, \dots,d_{m})$.

For any matrix $\mathbf{X} \in \mathfrak{F}^{n \times m}$, $d\mathbf{X}$ denotes the\emph{
matrix of differentials} $(dx_{ij})$. Finally, we define the \emph{measure} or volume element
$(d\mathbf{X})$ when $\mathbf{X} \in \mathfrak{F}^{m \times n}, \mathfrak{S}_{m}^{\beta}$,
$\mathfrak{D}_{m}^{\beta}$ or $\mathcal{V}_{m,n}^{\beta}$, see \citet{dggj:11}.

If $\mathbf{X} \in \mathfrak{F}^{m \times n}$ then $(d\mathbf{X})$ (the Lebesgue measure in
$\mathfrak{F}^{m \times n}$) denotes the exterior product of the $\beta mn$ functionally
independent variables
$$
  (d\mathbf{X}) = \bigwedge_{i = 1}^{m}\bigwedge_{j = 1}^{n}dx_{ij} \quad \mbox{ where }
    \quad dx_{ij} = \bigwedge_{k = 1}^{\beta}dx_{ij}^{(k)}.
$$

If $\mathbf{S} \in \mathfrak{S}_{m}^{\beta}$ (or $\mathbf{S} \in \mathfrak{T}_{U}^{\beta}(m)$
is a upper triangular matrix) then $(d\mathbf{S})$ (the Lebesgue measure in
$\mathfrak{S}_{m}^{\beta}$ or in $\mathfrak{T}_{U}^{\beta}(m)$) denotes the exterior product of
the exterior product of the $m(m-1)\beta/2 + m$ functionally independent variables,
$$
  (d\mathbf{S}) = \bigwedge_{i=1}^{m} ds_{ii}\bigwedge_{i < j}^{m}\bigwedge_{k = 1}^{\beta}
                      ds_{ij}^{(k)}.
$$
Observe, that for the Lebesgue measure $(d\mathbf{S})$ defined thus, it is required that
$\mathbf{S} \in \mathfrak{P}_{m}^{\beta}$, that is, $\mathbf{S}$ must be a non singular
Hermitian matrix (Hermitian definite positive matrix).

If $\mathbf{\Lambda} \in \mathfrak{D}_{m}^{\beta}$ then $(d\mathbf{\Lambda})$ (the Legesgue
measure in $\mathfrak{D}_{m}^{\beta}$) denotes the exterior product of the $\beta m$
functionally independent variables
$$
  (d\mathbf{\Lambda}) = \bigwedge_{i = 1}^{n}\bigwedge_{k = 1}^{\beta}d\lambda_{i}^{(k)}.
$$
If $\mathbf{H}_{1} \in \mathcal{V}_{m,n}^{\beta}$  is such that $\mathbf{H}_{1} =
(\mathbf{h}_{1}, \dots, \mathbf{h}_{m})$, where $\mathbf{h}_{j}$, $j = 1, \dots,m$ are their
columns, then
$$
  (\mathbf{H}^{*}_{1}d\mathbf{H}_{1}) = \bigwedge_{i=1}^{m} \bigwedge_{j =i+1}^{n}
  \mathbf{h}_{j}^{*}d\mathbf{h}_{i}.
$$
where the partitioned matrix $\mathbf{H} = (\mathbf{H}_{1}|\mathbf{H}_{2}) = (\mathbf{h}_{1},
\dots, \mathbf{h}_{m}|\mathbf{h}_{m+1}, \dots, \mathbf{h}_{n}) \in \mathfrak{U}^{\beta}(n)$,
with $\mathbf{H}_{2} = (\mathbf{h}_{m+1}, \dots, \mathbf{h}_{n})$. It can be proved that this
differential form does not depend on the choice of the $\mathbf{H}_{2}$ matrix. When $n = 1$;
$\mathcal{V}^{\beta}_{m,1}$ defines the unit sphere in $\mathfrak{F}^{m}$. This is, of course,
an $(m-1)\beta$- dimensional surface in $\mathfrak{F}^{m}$. When $n = m$ and denoting
$\mathbf{H}_{1}$ by $\mathbf{H}$, $(\mathbf{H}d\mathbf{H}^{*})$ is termed the \emph{Haar
measure} on $\mathfrak{U}^{\beta}(m)$.

The multivariate \emph{Gamma function} for the space $\mathfrak{S}_{m}^{\beta}$ denotes as
$\Gamma^{\beta}_{m}[a]$, is defined by
\begin{eqnarray*}
  \Gamma_{m}^{\beta}[a] &=& \displaystyle\int_{\mathbf{A} \in \mathfrak{P}_{m}^{\beta}}
  \etr\{-\mathbf{A}\} |\mathbf{A}|^{a-(m-1)\beta/2 - 1}(d\mathbf{A}) \\
    &=& \pi^{m(m-1)\beta/4}\displaystyle\prod_{i=1}^{m} \Gamma[a-(i-1)\beta/2],
\end{eqnarray*}
where $\etr(\cdot) = \exp(\tr(\cdot))$, $|\cdot|$  denotes the determinant and $\re(a)
> (m-1)\beta/2$, see \citet{gr:87}. This can be obtained as a particular case of the
\emph{generalised gamma function of weight $\kappa$} for the space $\mathfrak{S}^{\beta}_{m}$
with $\kappa = (k_{1}, k_{2}, \dots, k_{m})$, $k_{1}\geq k_{2}\geq \cdots \geq k_{m} \geq 0$,
$k_{1}, k_{2},\dots, k_{m}$ are nonnegative integers, taking $\kappa =(0,0,\dots,0)$ and which
for $\re(a) \geq (m-1)\beta/2 - k_{m}$ is defined by, see \citet{gr:87},
\begin{eqnarray}\label{int1}
  \Gamma_{m}^{\beta}[a,\kappa] &=& \displaystyle\int_{\mathbf{A} \in \mathfrak{P}_{m}^{\beta}}
  \etr\{-\mathbf{A}\} |\mathbf{A}|^{a-(m-1)\beta/2 - 1} q_{\kappa}(\mathbf{A}) (d\mathbf{A}) \\
&=& \pi^{m(m-1)\beta/4}\displaystyle\prod_{i=1}^{m} \Gamma[a + k_{i}
    -(i-1)\beta/2]\nonumber\\ \label{gammagen1}
&=& [a]_{\kappa}^{\beta} \Gamma_{m}^{\beta}[a],
\end{eqnarray}
where for $\mathbf{A} \in \mathfrak{S}_{m}^{\beta}$
\begin{equation}\label{hwv}
    q_{\kappa}(\mathbf{A}) = |\mathbf{A}_{m}|^{k_{m}}\prod_{i = 1}^{m-1}|\mathbf{A}_{i}|^{k_{i}-k_{i+1}}
\end{equation}
with $\mathbf{A}_{p} = (a_{rs})$, $r,s = 1, 2, \dots, p$, $p = 1,2, \dots, m$ is termed the
\emph{highest weight vector}, see \citet{gr:87}. In other branches of mathematics the
\textit{highest weight vector} $q_{\kappa}(\mathbf{A})$ is also termed the \emph{generalised
power} of $\mathbf{A}$ and is denoted as $\Delta_{\kappa}(\mathbf{A})$, see \citet{fk:94} and
\citet{hl:01}.

Some additional properties of $q_{\kappa}(\mathbf{A})$, which are immediate consequences of the
definition of $q_{\kappa}(\mathbf{A})$ and the following property 1, are:
\begin{enumerate}
  \item Let $\mathbf{A} = \mathbf{L}^{*}\mathbf{DL}$ be the L'DL decomposition of $\mathbf{A} \in \mathfrak{P}_{m}^{\beta}$,
        where $\mathbf{L} \in \mathfrak{T}_{U}^{\beta}(m)$ with $l_{ii} = 1$, $i = 1, 2, \ldots ,m$ and
        $\mathbf{D} = \diag(\lambda_{1}, \dots, \lambda_{m})$, $\lambda_{i} \geq 0$, $i = 1, 2, \ldots
        ,m$. Then
        \begin{equation}\label{qk1}
          q_{\kappa}(\mathbf{A}) = \prod_{i=1}^{m} \lambda_{i}^{k_{i}}.
        \end{equation}
      \item
      \begin{equation}\label{qk2}
        q_{\kappa}(\mathbf{A}^{-1}) =  q_{-\kappa^{*}}^{*}(\mathbf{A}),
      \end{equation}
      where $\kappa^{*}=(k_{m}, k_{m-1}, \dots,k_{1})$, $-\kappa^{*}=(-k_{m}, -k_{m-1},
      \dots,-k_{1})$,
      \begin{equation}\label{hhwv}
         q_{\kappa}^{*}(\mathbf{A}) = |\mathbf{A}_{m}|^{k_{m}}\prod_{i = 1}^{m-1}|\mathbf{A}_{i}|^{k_{i}-k_{i+1}}
      \end{equation}
      and
      \begin{equation}\label{qqk1}
        q_{\kappa}^{*}(\mathbf{A}) = \prod_{i=1}^{m} \lambda_{i}^{k_{m-i+1}},
      \end{equation}
      see \citet[pp. 126-127 and Proposition VII.1.5]{fk:94}.

      Alternatively, let $\mathbf{A} = \mathbf{T}^{*}\mathbf{T}$ the Cholesky's decomposition of
  matrix $\mathbf{A} \in \mathfrak{P}_{m}^{\beta}$, with $\mathbf{T}=(t_{ij}) \in
  \mathfrak{T}_{U}^{\beta}(m)$, then $\lambda_{i} = t_{ii}^{2}$, $t_{ii} \geq 0$, $i = 1, 2,
  \ldots ,m$. See \citet[p. 931, first paragraph]{hl:01}, \citet[p. 390, lines -11 to
  -16]{hlz:05} and \citet[p.5, lines 1-6]{k:14}.
  \item if $\kappa = (p, \dots, p)$, then
    \begin{equation}\label{qk3}
        q_{\kappa}(\mathbf{A}) = |\mathbf{A}|^{p},
    \end{equation}
    in particular if $p=0$, then $q_{\kappa}(\mathbf{A}) = 1$.
  \item if $\tau = (t_{1}, t_{2}, \dots, t_{m})$, $t_{1}\geq t_{2}\geq \cdots \geq t_{m} \geq
  0$, then
    \begin{equation}\label{qk41}
        q_{\kappa+\tau}(\mathbf{A}) = q_{\kappa}(\mathbf{A})q_{\tau}(\mathbf{A}),
    \end{equation}
    in particular if $\tau = (p,p, \dots, p)$,  then
    \begin{equation}\label{qk42}
        q_{\kappa+\tau}(\mathbf{A}) \equiv q_{\kappa+p}(\mathbf{A}) = |\mathbf{A}|^{p} q_{\kappa}(\mathbf{A}).
    \end{equation}
    \item Finally, for $\mathbf{B} \in \mathfrak{T}_{U}^{\beta}(m)$  in such a manner that
$\mathbf{C} = \mathbf{B}^{*}\mathbf{B} \in \mathfrak{P}_{m}^{\beta}$,
    \begin{equation}\label{qk5}
        q_{\kappa}(\mathbf{B}^{*}\mathbf{AB}) = q_{\kappa}(\mathbf{C})q_{\kappa}(\mathbf{A}),
    \end{equation}
see \citet[p. 776, eq. (2.1)]{hlz:08}.
\end{enumerate}
\begin{remark}
Formally the determinant is defined in terms of function $q_{\kappa}(\mathbf{A})$, more
precisely in terms of the function $q(\mathbf{A})$ which is defined as:
$$
  q(\mathbf{A}) = |\mathbf{A}|^{\eta}
$$
where
$$
  \eta = \left\{
        \begin{array}{ll}
          1 & \hbox{if } \mathfrak{F} = \Re \hbox{ or } \mathfrak{C} \\
          1/2 & \hbox{if } \mathfrak{F} = \mathfrak{H}\\
          1/4 & \hbox{if } \mathfrak{F} = \mathfrak{C}.
        \end{array}
      \right.
$$
Alternatively, by (\ref{qk1}) the determinant of a matrix $\mathbf{A}$ is defined in terms of
their eigenvalues, see \citet{gr:87}, \citet{fk:94} and \citet{f:05}.
\end{remark}

\begin{remark}
Let $\mathcal{P}(\mathfrak{S}_{m}^{\beta})$ denote the algebra of all polynomial functions on
$\mathfrak{S}_{m}^{\beta}$, and $\mathcal{P}_{k}(\mathfrak{S}_{m}^{\beta})$ the subspace of
homogeneous polynomials of degree $k$ and let $\mathcal{P}^{\kappa}(\mathfrak{S}_{m}^{\beta})$
be an irreducible subspace of $\mathcal{P}(\mathfrak{S}_{m}^{\beta})$ such that
$$
  \mathcal{P}_{k}(\mathfrak{S}_{m}^{\beta}) = \sum_{\kappa}\bigoplus
  \mathcal{P}^{\kappa}(\mathfrak{S}_{m}^{\beta}).
$$
Note that $q_{\kappa}$ is a homogeneous polynomial of degree $k$, moreover $q_{\kappa} \in
\mathcal{P}^{\kappa}(\mathfrak{S}_{m}^{\beta})$, see \citet{gr:87}.
\end{remark}
In (\ref{gammagen1}), $[a]_{\kappa}^{\beta}$ denotes the generalised Pochhammer symbol of
weight $\kappa$, defined as
$$
  [a]_{\kappa}^{\beta} = \prod_{i = 1}^{m}(a-(i-1)\beta/2)_{k_{i}}
    = \frac{\Gamma_{m}^{\beta}[a,\kappa]}{\Gamma_{m}^{\beta}[a]},
$$
where $\re(a) > (m-1)\beta/2 - k_{m}$ and $(a)_{i} = a (a+1)\cdots(a+i-1),$ is the standard
Pochhammer symbol.

A variant of the generalised gamma function of weight $\kappa$ is obtained from
\citet{k:66} and is defined as%
\begin{eqnarray}\label{int2}
  \Gamma_{m}^{\beta}[a,-\kappa] &=& \displaystyle\int_{\mathbf{A} \in \mathfrak{P}_{m}^{\beta}}
    \etr\{-\mathbf{A}\} |\mathbf{A}|^{a-(m-1)\beta/2 - 1} q_{\kappa}(\mathbf{A}^{-1})
    (d\mathbf{A}) \\
&=& \pi^{m(m-1)\beta/4}\displaystyle\prod_{i=1}^{m} \Gamma[a - k_{i}
    -(i-1)\beta/2] \nonumber \\ \label{gammagen2}
&=& \displaystyle\frac{(-1)^{k} \Gamma_{m}^{\beta}[a]}{[-a +(m-1)\beta/2
    +1]_{\kappa}^{\beta}} ,
\end{eqnarray}
where $\re(a) > (m-1)\beta/2 + k_{1}$.

Now, we show three Jacobians in terms of the $\beta$ parameter, which are proposed as
extensions of real, complex or quaternion cases, see \citet{dggj:11}.

\begin{lemma}\label{lemhlt} Let $\mathbf{X}$ and $\mathbf{Y} \in
\mathfrak{P}_{m}^{\beta}$ matrices of functionally independent variables, and let $\mathbf{Y} =
\mathbf{AXA^{*}} + \mathbf{C}$, where $\mathbf{A} \in {\mathcal L}_{m,m}^{\beta} $ and
$\mathbf{C} \in \mathfrak{P}_{m}^{\beta}$ are matrices of constants. Then
\begin{equation}\label{hlt}
    (d\mathbf{Y}) = |\mathbf{A}^{*}\mathbf{A}|^{\beta(m-1)/2+1} (d\mathbf{X}).
\end{equation}
\end{lemma}

\begin{lemma}\label{lemi}
Let $\mathbf{S} \in \mathfrak{P}_{m}^{\beta}.$ Then ignoring the sign, if $\mathbf{Y} =
\mathbf{S}^{-1}$
\begin{equation}\label{i}
    (d\mathbf{Y}) = |\mathbf{S}|^{-\beta(m - 1) - 2}(d\mathbf{S}).
\end{equation}
\end{lemma}
We end this section with a some general results, which are useful in a variety of situations,
which enable us to transform the density function of a matrix $\mathbf{X} \in
\mathfrak{P}_{m}^{\beta}$ to the density function of its eigenvalues, see \citet{dggj:09b}.

\begin{lemma}\label{eig}
Let $\mathbf{X} \in \mathfrak{P}_{m}^{\beta}$ be a random matrix with density function
$f_{_{\mathbf{X}}}(\mathbf{X})$. Then the joint density function of the eigenvalues
$\lambda_{1}, \dots, \lambda_{m}$ of $\mathbf{X}$ is
\begin{equation}\label{dfeig}
    \frac{\pi^{m^{2}\beta/2+ \varrho}}{\Gamma_{m}^{\beta}[m\beta/2]} \prod_{i<
    j}^{m}(\lambda_{i} -\lambda_{j})^{\beta}\int_{\mathbf{H} \in
    \mathfrak{U}^{\beta}(m)}f_{_{\mathbf{X}}}(\mathbf{HLH}^{*})(d\mathbf{H})
\end{equation}
where $\mathbf{L} = \diag(\lambda_{1}, \dots, \lambda_{m})$, $\lambda_{1}> \cdots > \lambda_{m}
> 0$, $(d\mathbf{H})$ is the normalised Haar measure and
$$
  \varrho = \left\{
             \begin{array}{rl}
               0, & \beta = 1; \\
               -m, & \beta = 2; \\
               -2m, & \beta = 4; \\
               -4m, & \beta = 8.
             \end{array}
           \right.
$$
\end{lemma}

\section{Riesz distributions}\label{sec3}

Alternatively to find the Riesz distribution from the Riesz measure, \citet{hl:01}, in this
section we establish the two versions of the Riesz distribution via the next result. For this
purpose, we utilise the complexification $\mathfrak{S}_{m}^{\beta, \mathfrak{C}} =
\mathfrak{S}_{m}^{\beta} + i \mathfrak{S}_{m}^{\beta}$ of $\mathfrak{S}_{m}^{\beta}$. That is,
$\mathfrak{S}_{m}^{\beta, \mathfrak{C}}$ consist of all matrices $\mathbf{Z} \in
(\mathfrak{F^{\mathfrak{C}}})^{m \times m}$ of the form $\mathbf{Z} = \mathbf{X} +
i\mathbf{Y}$, with $\mathbf{X}, \mathbf{Y} \in \mathfrak{S}_{m}^{\beta}$. We refer to
$\mathbf{X} = \re(\mathbf{Z})$ and $\mathbf{Y} = \im(\mathbf{Z})$ as the \emph{real and
imaginary parts} of $\mathbf{Z}$, respectively. The \emph{generalised right half-plane}
$\mathbf{\Phi}_{m}^{\beta} = \mathfrak{P}_{m}^{\beta} + i \mathfrak{S}_{m}^{\beta}$ in
$\mathfrak{S}_{m}^{\beta,\mathfrak{C}}$ consists of all $\mathbf{Z} \in
\mathfrak{S}_{m}^{\beta,\mathfrak{C}}$ such that $\re(\mathbf{Z}) \in
\mathfrak{P}_{m}^{\beta}$, see \cite[p. 801]{gr:87}.

\begin{lemma}\label{lemmaC}
Let $\mathbf{\Sigma} \in \mathbf{\Phi}_{m}^{\beta}$ and  $\kappa = (k_{1}, k_{2}, \dots,
k_{m})$, $k_{1}\geq k_{2}\geq \cdots \geq k_{m} \geq 0$, $k_{1}, k_{2},\dots, k_{m}$ are
nonnegative integers. Then
\begin{enumerate}
  \item for $\re(a) \geq (m-1)\beta/2 - k_{m}$,
  $$
    \displaystyle\int_{\mathbf{A} \in \mathfrak{P}_{m}^{\beta}} \etr\{-\beta \mathbf{\Sigma}^{-1}\mathbf{A}\}
  |\mathbf{A}|^{a-(m-1)\beta/2 - 1} q_{\kappa}(\mathbf{A}) (d\mathbf{A}) \hspace{3cm}
  $$
  \begin{equation}\label{intR1}
   \hspace{4cm}= \displaystyle\frac{\Gamma_{m}^{\beta}[a,\kappa] |\mathbf{\Sigma}|^{a}q_{\kappa}(\mathbf{\Sigma})}
   {\beta^{am+\sum_{i = 1}^{m}k_{i}}}
  \end{equation}
  \item for $\re(a) > (m-1)\beta/2 + k_{1}$,
  $$
    \displaystyle\int_{\mathbf{A} \in \mathfrak{P}_{m}^{\beta}} \etr\{-\beta \mathbf{\Sigma}^{-1}\mathbf{A}\}
  |\mathbf{A}|^{a-(m-1)\beta/2 - 1} q_{\kappa}(\mathbf{A}^{-1}) (d\mathbf{A})\hspace{3cm}
  $$
  \begin{equation}\label{intR2}
   \hspace{4.2cm}= \frac{\Gamma_{m}^{\beta}[a,-\kappa]
   |\mathbf{\Sigma}|^{a}q_{\kappa}(\mathbf{\Sigma}^{-1})}{\beta^{am-\sum_{i = 1}^{m}k_{i}}}
  \end{equation}
\end{enumerate}
\end{lemma}
\begin{proof}
Let $\mathbf{\Sigma} \in \mathbf{\Phi}_{m}^{\beta}$ and in integrals (\ref{intR1}) and
(\ref{intR2}) make the change of variable $\mathbf{A} = \left(\mathbf{\Sigma}^{1/2}\right)^{*}
\mathbf{W \Sigma}^{1/2}/\beta$, where $\mathbf{\Sigma}^{1/2} \in \mathfrak{T}_{U}^{\beta}(m)$
denotes the Cholesky decomposition of $\mathbf{\Sigma}$, such that
$\left(\mathbf{\Sigma}^{1/2}\right)^{*} \mathbf{\Sigma}^{1/2} = \mathbf{\Sigma}$, see
\citet[Subsection 1.7, p. 788]{gr:87} and \citet[Chapter VI, p. 100]{fk:94}. By Lemma
\ref{lemhlt},
$$
  (d\mathbf{A}) = \beta^{-m(m-1)\beta/2-m} |\mathbf{\Sigma}|^{\beta(m-1)/2+1}
(d\mathbf{W}),
$$
then integrals becomes
\begin{enumerate}
  \item
  $$
  \beta^{-ma}\int_{\mathbf{W} \in \mathfrak{P}_{m}^{\beta}} \etr\{-\mathbf{W}\}
  |\mathbf{W}|^{a-(m-1)\beta/2 - 1} q_{\kappa}\left(\frac{1}{\beta}\left(\mathbf{\Sigma}^{1/2}\right)^{*} \mathbf{W \Sigma}^{1/2}\right)
  |\mathbf{\Sigma}|^{a} (d\mathbf{W}).
  $$
  \item and
  $$
  \beta^{-ma}\int_{\mathbf{W} \in \mathfrak{P}_{m}^{\beta}} \etr\{-\mathbf{W}\}
  |\mathbf{W}|^{a-(m-1)\beta/2 - 1} q_{\kappa}(\beta\mathbf{\Sigma}^{-1/2} \mathbf{W}^{-1}\left(\mathbf{\Sigma}^{-1/2}\right)^{*})
  |\mathbf{\Sigma}|^{a} (d\mathbf{W}).
  $$
\end{enumerate}
The desired results now follows from (\ref{int1}) and (\ref{qk5}), and (\ref{int2}) ,
respectively.
\end{proof}

Hence, as consequence of Lemma \ref{lemmaC}, now we can propose the two definitions of Riesz
distribution.

\begin{definition}\label{defnRd}
Let $\mathbf{\Sigma} \in \mathbf{\Phi}_{m}^{\beta}$ and  $\kappa = (k_{1}, k_{2}, \dots,
k_{m})$, $k_{1}\geq k_{2}\geq \cdots \geq k_{m} \geq 0$, $k_{1}, k_{2},\dots, k_{m}$ are
nonnegative integers.
\begin{enumerate}
  \item Then it said that $\mathbf{X}$ has a Riesz distribution of type I if its density function is
  \begin{equation}\label{dfR1}
    \frac{\beta^{am+\sum_{i = 1}^{m}k_{i}}}{\Gamma_{m}^{\beta}[a,\kappa] |\mathbf{\Sigma}|^{a}q_{\kappa}(\mathbf{\Sigma})}
    \etr\{-\beta \mathbf{\Sigma}^{-1}\mathbf{X}\}|\mathbf{X}|^{a-(m-1)\beta/2 - 1}
    q_{\kappa}(\mathbf{X})(d\mathbf{X})
  \end{equation}
  for $\mathbf{X} \in \mathfrak{P}_{m}^{\beta}$ and $\re(a) > (m-1)\beta/2 - k_{m}$.
  \item Then it said that $\mathbf{X}$ has a Riesz distribution of type II if its density function is
  \begin{equation}\label{dfR2}
     \frac{\beta^{am-\sum_{i = 1}^{m}k_{i}}}{\Gamma_{m}^{\beta}[a,-\kappa]
   |\mathbf{\Sigma}|^{a}q_{\kappa}(\mathbf{\Sigma}^{-1})}\etr\{-\beta \mathbf{\Sigma}^{-1}\mathbf{X}\}
  |\mathbf{X}|^{a-(m-1)\beta/2 - 1} q_{\kappa}(\mathbf{X}^{-1}) (d\mathbf{X})
  \end{equation}
  for $\mathbf{X} \in \mathfrak{P}_{m}^{\beta}$ and $\re(a) > (m-1)\beta/2 + k_{1}$.
\end{enumerate}
\end{definition}

Note that, the matrix multivariate gamma distribution is a particular example of the Riesz
distribution. Furthermore, if $\kappa = (0,0, \dots,0)$ in two densities in Definition
\ref{defnRd} the matrix multivariate gamma distribution is obtained. In addition observe that,
in this last case $\mathbf{\Sigma} = 2\mathbf{\Sigma}$ and $a = \beta n/2$, the Wishart
distribution is gotten.

In next result is established their characteristic functions of Riesz distributions.

\begin{theorem}\label{theoRd}
Let $\mathbf{\Sigma} \in \mathbf{\Phi}_{m}^{\beta}$ and  $\kappa = (k_{1}, k_{2}, \dots,
k_{m})$, $k_{1}\geq k_{2}\geq \cdots \geq k_{m} \geq 0$, $k_{1}, k_{2},\dots, k_{m}$ are
nonnegative integers.
\begin{enumerate}
  \item Then if $\mathbf{X}$ has a Riesz distribution of type I its characteristic function is
  \begin{equation}\label{chfR1}
    |\mathbf{I}_{m}-i\mathbf{\Sigma T}/\beta|^{-a}  q_{\kappa}\left(\left(\mathbf{I}_{m}-i\mathbf{\Sigma}^{1/2}
    \mathbf{T}\left(\mathbf{\Sigma}^{1/2}\right)^{*}/\beta\right)^{-1}\right)
  \end{equation}
  for $\re(a) \geq (m-1)\beta/2 - k_{m}$.
  \item Then if $\mathbf{X}$ has a Riesz distribution of type II its characteristic function is
  \begin{equation}\label{chfR2}
     |\mathbf{I}_{m}-i\mathbf{\Sigma T}/\beta|^{-a}
    q_{\kappa}\left(\left(\mathbf{I}_{m}-i\mathbf{\Sigma}^{1/2}\mathbf{T}\left(\mathbf{\Sigma}^{1/2}\right)^{*}/\beta\right)\right)
  \end{equation}
  for $\re(a) > (m-1)\beta/2 + k_{1}$.
\end{enumerate}
\end{theorem}
\begin{proof}
1. The characteristic function is
$$
    \frac{\beta^{am+\sum_{i = 1}^{m}k_{i}}}{\Gamma_{m}^{\beta}[a,\kappa] |\mathbf{\Sigma}|^{a}q_{\kappa}(\mathbf{\Sigma})}
    \int_{\mathbf{X} \in \mathfrak{P}_{m}^{\beta}}\etr\{-\beta(\mathbf{\Sigma}^{-1}-i\mathbf{T}/\beta)
    \mathbf{X}\}|\mathbf{X}|^{a-(m-1)\beta/2 - 1} q_{\kappa}(\mathbf{X})(d\mathbf{X}).
$$
From (\ref{intR1}) the characteristic function is
$$
   \frac{|\mathbf{\Sigma}^{-1}-i\mathbf{T}/\beta|^{-a}q_{\kappa}((\mathbf{\Sigma}^{-1}-i\mathbf{T}/\beta)^{-1})}
     {|\mathbf{\Sigma}|^{a}q_{\kappa}(\mathbf{\Sigma})} = \frac{|\mathbf{I}_{m}-i\mathbf{T\Sigma}/\beta|^{-a}
     q_{\kappa}((\mathbf{\Sigma}^{-1}-i\mathbf{T}/\beta)^{-1})}{q_{\kappa}(\mathbf{\Sigma})}
$$
The final expression is obtained from (\ref{qk5}). The results stated in 2 is obtained in
analogously.
\end{proof}

Finally, it is found the joint distributions of the eigenvalues for random matrices Riesz type
I and II.

\begin{theorem}\label{teoE}
Let $\mathbf{\Sigma} \in \mathbf{\Phi}_{m}^{\beta}$ and  $\kappa = (k_{1}, k_{2}, \dots,
k_{m})$, $k_{1}\geq k_{2}\geq \cdots \geq k_{m} \geq 0$, $k_{1}, k_{2},\dots, k_{m}$ are
nonnegative integers.
\begin{enumerate}
  \item Let $\lambda_{1}, \dots,\lambda_{m}$,  $\lambda_{1}> \cdots >\lambda_{m} > 0$ be
  the eigenvalues of $\mathbf{X}$. Then if $\mathbf{X}$ has a Riesz distribution of type I,
  the joint density  of $\lambda_{1}, \dots, \lambda_{m}$ is
  $$
     \frac{\beta^{am+\sum_{i = 1}^{m}k_{i}} \ \pi^{m^{2}\beta/2+ \varrho}}{\Gamma_{m}^{\beta}[m\beta/2]\Gamma_{m}^{\beta}[a,\kappa]
     |\mathbf{\Sigma}|^{a}q_{\kappa}(\mathbf{\Sigma})} \prod_{i< j}^{m}(\lambda_{i} -\lambda_{j})^{\beta}
     \hspace{4cm}
  $$
  \begin{equation}\label{dER1}\hspace{4cm}
     \times \ \prod_{i=1}^{m}\lambda_{i}^{a-(m-1)\beta/2 - 1} \G_{\kappa}^{(m),\beta}(-\beta\mathbf{\Sigma}^{-1},\mathbf{L}).
  \end{equation}
  where $\mathbf{L} = \diag(\lambda_{1}, \dots,\lambda_{m})$, $\re(a) \geq (m-1)\beta/2 -
  k_{m}$ and
  $$
    \G_{\kappa}^{(m),\beta}(\mathbf{A},\mathbf{B}) =\int_{\mathbf{H} \in \mathfrak{U}^{\beta}(m)}
    \etr\{\mathbf{A}\mathbf{HBH}^{*}\} q_{\kappa}(\mathbf{HBH}^{*})(d\mathbf{H}),
  $$
  \item Let $\delta_{1}, \dots,\delta_{m}$,  $\delta_{1}> \cdots >\delta_{m} > 0$ be
  the eigenvalues of $\mathbf{X}$. Then if $\mathbf{X}$ has a Riesz distribution of type II,
  the joint density of their eigenvalues is
  $$
     \frac{\beta^{am-\sum_{i = 1}^{m}k_{i}} \ \pi^{m^{2}\beta/2+ \varrho}}
     {\Gamma_{m}^{\beta}[m\beta/2]\Gamma_{m}^{\beta}[a,\kappa]
     |\mathbf{\Sigma}|^{a}q_{\kappa}(\mathbf{\Sigma}^{-1})} \prod_{i< j}^{m}(\delta_{i} -\delta_{j})^{\beta}
     \hspace{4cm}
  $$
  \begin{equation}\label{dER2}\hspace{4cm}
     \times \ \prod_{i=1}^{m}\delta_{i}^{a-(m-1)\beta/2 - 1} \Q_{\kappa}^{(m),\beta}\left(-\beta\mathbf{\Sigma}^{-1},
     \mathbf{D}\right).
  \end{equation}
  where $\mathbf{D} = \diag(\delta_{1}, \dots,\delta_{m})$, $\re(a) > (m-1)\beta/2 + k_{1}$ and
  $$
    \Q_{\kappa}^{(m),\beta}(\mathbf{A},\mathbf{B}) =\int_{\mathbf{H} \in \mathfrak{U}^{\beta}(m)}
    \etr\{\mathbf{A}\mathbf{HBH}^{*}\} q_{\kappa}(\mathbf{HB}^{-1}\mathbf{H}^{*})(d\mathbf{H}).
  $$
\end{enumerate}
And $\varrho$ is defined in Lemma 3.
\end{theorem}
\begin{proof} 1. From Lemma 3
$$
    \frac{\beta^{am+\sum_{i = 1}^{m}k_{i}} \ \pi^{m^{2}\beta/2+ \varrho}}{\Gamma_{m}^{\beta}[m\beta/2]\Gamma_{m}^{\beta}[a,\kappa]
    |\mathbf{\Sigma}|^{a}q_{\kappa}(\mathbf{\Sigma})} \prod_{i<
    j}^{m}(\lambda_{i} -\lambda_{j})^{\beta} \hspace{4cm}
$$
$$\hspace{2cm}
    \int_{\mathbf{H} \in \mathfrak{U}^{\beta}(m)}\etr\{-\beta \mathbf{\Sigma}^{-1}\mathbf{HLH}^{*}\}
    |\mathbf{HLH}^{*}|^{a-(m-1)\beta/2 - 1} q_{\kappa}(\mathbf{HLH}^{*})(d\mathbf{H}).
$$
Therefore,
$$
    \frac{\beta^{am+\sum_{i = 1}^{m}k_{i}} \ \pi^{m^{2}\beta/2+ \varrho}}{\Gamma_{m}^{\beta}[m\beta/2]\Gamma_{m}^{\beta}[a,\kappa]
    |\mathbf{\Sigma}|^{a}q_{\kappa}(\mathbf{\Sigma})} \prod_{i<
    j}^{m}(\lambda_{i} -\lambda_{j})^{\beta} \prod_{i=1}^{m}\lambda_{i}^{a-(m-1)\beta/2 - 1}
    \hspace{2cm}
$$
$$\hspace{5.5cm}
    \int_{\mathbf{H} \in \mathfrak{U}^{\beta}(m)}\etr\{-\beta \mathbf{\Sigma}^{-1}\mathbf{HLH}^{*}\}
     q_{\kappa}(\mathbf{HLH}^{*})(d\mathbf{H}).
$$
2. is proved similarly.
\end{proof}

Again, observe that when $\kappa = (0, \dots,0)$, the corresponding joint densities of the
eigenvalues of the matrix multivariate gamma and Wishart distributions are obtained as
particular cases of (\ref{dER1}) and (\ref{dER2}), indistinctly.

\section{Conclusions}

Potentially, the Riesz distribution can replace the Wishart distribution throughout the
multivariate statistical analysis, however still needs to find the sampling distribution
associated Riesz distribution, i.e. the distribution of the matrix $\mathbf {X}$ such that the
random matrix $\mathbf{V} = \mathbf{X}^{*}\mathbf{X}$ has a distribution of Riesz and, thus
appears naturally the Riesz distribution instead of the Wishart distribution. This topic is
currently being studied by the author. But the Wishart distribution occurs naturally in other
fields of knowledge without being associated with the matrix multivariate normal distribution,
as in the case of the joint element densities of generalised Laguerre ensembles. Thus, if the
densities (\ref{dfR1}) and (\ref{dfR2}) are considered as the joint element densities of
generalised Laguerre ensembles, the distributions (\ref{dER1}) and (\ref{dER2}) are the
distributions of the Laguerre ensembles under a Riesz distributions, see \cite{dggj:13}.

Any topic in statistical literature, is usually first studied in the real case, later in the
complex case, later for the quaternion case and exceptionally for the octonion case. From the
results presented in this paper, the real, complex, quaternion and octonion cases are obtained
by simply replacing $\beta$ with $1,2,4$ or $8$, respectively.

Finally, note that the real dimension of real normed division algebras can be expressed as
powers of 2, $\beta = 2^{n}$ for $n = 0,1,2,3$. On the other hand, as it can be check in
\citet{k:84}, the results obtained in this work can be extended to the hypercomplex cases; that
is, for complex, bicomplex, biquaternion and bioctonion (or sedenionic) algebras, which of
course are not division algebras (except the complex algebra). Note, also, that hypercomplex
algebras are obtained by replacing the real numbers with complex numbers in the construction of
real normed division algebras. Thus, the results for hypercomplex algebras are obtained by
simply replacing $\beta$ with $2\beta$ in our results. Alternatively, following \citet{k:84},
we can conclude that, our results are true for `$2^{n}$-ions', $n = 0,1,2,3,4,5$, emphasising
that only for $n=0,1,2,3$ the corresponding algebras are real normed division algebras.

Here we prefer our notation used unlike the proposal in \citet{hl:01}, because we think that is
clearer to a statistician. Note, however, that any of the results obtained can be  written
easily as proposed by \citet{hl:01} in the context of Euclidean simple Jordan algebras, simply
applying (\ref{qk42}). For example in this case (\ref{chfR1}) can be written equivalently as:
$$
    q_{\kappa+a}\left(\left(\mathbf{I}_{m}-i\mathbf{\Sigma}^{1/2}
    \mathbf{T}\left(\mathbf{\Sigma}^{1/2}\right)^{*}/\beta\right)^{-1}\right).
$$




\begin{thebibliography}{99}
\bibitem[Baez(2002)]{b:02}
    J. C. Baez,
    The octonions,
    Bull. Amer. Math. Soc.
    39 (2002) 145--205.

\bibitem[Bhavsar(2000)]{b:00}
    C. D. Bhavsar,
    Asymptotic distributions of likelihood ratio criteria for two testing problems.
    Kybernetes 29(4)(2002) 510--517.

\bibitem[Casalis and Letac(1996)]{cl:96}
    M. Casalis, G. Letac,
    The Lukascs-Olkin-Rubin characterization of Wishart distributions on symmetric cones,
    Ann. Statist.  24 (1996) 768--786.

\bibitem[Constantine(1963)]{c:63}
   A. C. Constantine,
   Noncentral distribution problems in multivariate analysis.
   Ann. Math. Statist. 34 (1963) 1270--1285.

\bibitem[Davis(1980)]{da:80}
   A. W. Davis,
   Invariant polynomials with two matrix arguments, extending the zonal polynomials,
   In: Krishnaiah P R (ed.) Multivariate Analysis V.
   North-Holland Publishing Company, pp. 287--299, 1980.

\bibitem[D\'{\i}az-Garc\'{\i}a(2011)]{dg:11}
    J. A. D\'{\i}az-Garc\'{\i}a,
    Generalizations of some properties of invariant polynomials with matrix arguments,
    Appl. Math. (Warsaw) 38(4)(2011), 469-475.

\bibitem[D\'{\i}az-Garc\'{\i}a and  Guti\'errez-J\'aimez(2009)]{dggj:09b}
    J. A. D\'{\i}az-Garc\'{\i}a, R. Guti\'errez-J\'aimez,
    Special functions: Integral properties of Jack polynomials, hypergeometric
    functions and Invariant polynomials,
    \texttt{http://arxiv.org/abs/0909.1988}, 2009. Also submited.

\bibitem[D\'{\i}az-Garc\'{\i}a and  Guti\'errez-J\'aimez(2011)]{dggj:11}
    J. A. D\'{\i}az-Garc\'{\i}a,  and  R. Guti\'errez-J\'aimez,
    On Wishart distribution: some extensions,
    Linear Algebra Appl. 435 (2011) 1296-1310.

\bibitem[D\'{\i}az-Garc\'{\i}a and  Guti\'errez-J\'aimez(2012a)]{dggj:12a}
    J. A. D\'{\i}az-Garc\'{\i}a,  and  R. Guti\'errez-J\'aimez,
     Matricvariate and matrix multivariate T distributions and associated distributions,
     Metrika, 75(7)(2012a) 963-976.

\bibitem[D\'{\i}az-Garc\'{\i}a and  Guti\'errez-J\'aimez(2012b)]{dggj:12b}
    J. A. D\'{\i}az-Garc\'{\i}a,  and  R. Guti\'errez-J\'aimez,
    An identity of Jack polynomials,
    J. Iranian Statist. Soc. 11(1)(2012b) 87-92.

\bibitem[D\'{\i}az-Garc\'{\i}a and  Guti\'errez-J\'aimez(2013)]{dggj:13}
     J. A. D\'{\i}az-Garc\'{\i}a, R. Guti\'errez-J\'aimez,
     Spherical ensembles.
     Linear Algebra Appl. 438 (2013)  3174-3201.

\bibitem[Faraut and Kor\'anyi(1994)]{fk:94}
    J. Faraut, A. Kor\'anyi,
    Analysis on symmetric cones,
    Oxford Mathematical Monographs,
    Clarendon Press, Oxford, 1994.

\bibitem[Forrester (2009)]{f:05}
    P. J. Forrester,
    Log-gases and random matrices,
    To appear. (Available in: \verb"http://www.ms.unimelb.edu.au/~matpjf/matpjf.html".

\bibitem[Goodman (1963)]{g:63}
   N. R. Goodman,
   Statistical analysis based on a certain multivariate complex Gaussian distribution (an
   introdiction).
   Ann. Math. Statist.
   34(1)(1963) 152--177.

\bibitem[Gross and Richards(1987)]{gr:87}
    K. I. Gross, D. ST. P. Richards,
    Special functions of matrix argument I: Algebraic induction zonal polynomials
    and hypergeometric  functions,
    Trans. Amer. Math. Soc.
    301(2) (1987) 475--501.

\bibitem[Hassairi and Lajmi(2001)]{hl:01}
   A. Hassairi, S. Lajmi,
   Riesz exponential families on symmetric cones,
   J. Theoret. Probab. 14 (2001) 927-–948.

\bibitem[Hassairi \textit{et al.}(2005)]{hlz:05}
   A. Hassairi, S. Lajmi, R. Zine,
   Beta-Riesz distributions on symmetric cones,
   J. Statist. Plann. Inf. 133 (2005) 387–-404.

\bibitem[Hassairi \textit{et al.}(2008)]{hlz:08}
   A. Hassairi, S. Lajmi, R. Zine,
   A chacterization of the Riesz probability distribution,
   J. Theoret. Probab. 21 (2008) 773-–790.

\bibitem[Herz(1955)]{h:55}
   C. S. Herz,
   Bessel functions of matrix argument.
   Ann. of Math.
   61(3)(1955) 474-523.

\bibitem[James(1961)]{j:61}
    A. T. James,
    Zonal polynomials of the real positive definite symmetric matrices,
    Ann. Math. 35 (1961) 456--469.

\bibitem[James(1964)]{j:64}
    A. T. James,
    Distribution of matrix variate and latent roots derived from normal samples,
    Ann. Math. Statist.
    35 (1964) 475--501.

\bibitem[Kabe(1984)]{k:84}
    D. G. Kabe,
    Classical statistical analysis based on a certain hypercomplex multivariate
    normal distribution,
    \textit{Metrika}
    \textbf{31}(1984) 63--76.

\bibitem[Khatri(1966)]{k:66}
    C. G. Khatri,
    On certain distribution problems based on positive definite quadratic
   functions in normal vector,
    Ann. Math. Statist.
    37 (1966) 468--479.

\bibitem[Ko{\l}odziejek(2014)]{k:14}
   B. Ko{\l}odziejek
   The Lukacs-Olkin-Rubin theorem on symmetric cones without invariance of the ``Quotient".
   J. Theoret. Probab. (2014), DOI 10.1007/s10959-014-0587-3.

\bibitem[Li and Xue (2009)]{lx:09}
    F. Li, Y. Xue,
    Zonal polynomials and hypergeometric functions of quaternion matrix argument.
    Comm. Statist. Theory Methods
    38(8) (2009) 1184-1206.

\bibitem[Massam(1994)]{m:94}
    H. Massam,
    An exact decomposition theorem and unified view of some related distributions for a class
    of exponential transformation models on symmetric cones,
    Ann. Statist. 22(1)(1994) 369--394.

\bibitem[Micheas \emph{et al.} (2006)]{mdm:06}
    A. C. Micheas, D. K. Dey, K. V. Mardia
    Complex elliptical distribution with application to shape theory.
    J. Statist. Plann. Infer.  136 (2006) 2961-2982.

\bibitem[Muirhead(1982)]{m:82}
    R. J.Muirhead,
    Aspects of Multivariate Statistical Theory,
    John Wiley \& Sons, New York, 1982.

\bibitem[Riesz(1949)]{r:49}
    M. Riesz,
    L'int\'egrale de Riemann-Liouville et le probl\`{e}me de Cauchy.
    Acta Math. 81 (1949) 1--223.

\bibitem[Sawyer(1997)]{s:97}
    P. Sawyer,
    Spherical Functions on Symmetric Cones,
    Trans. Amer. Math. Soc.
    349 (1997) 3569--3584.

\bibitem[Wooding (1956)]{w:56}
    R. A. Wooding
    The multivariate distribution of complex normal variables.
    Biometrika
    43(1) (1956) 212--215.
\end{thebibliography}
\end{document}